\documentclass[10pt,reqno]{article}
\usepackage{amsmath,amsthm,amsfonts,amssymb}
\usepackage{epsfig}
\usepackage{color}
\usepackage{latexsym}
\usepackage{supertabular}
\usepackage{graphicx}
\usepackage{a4wide}
\numberwithin{equation}{section}

\def\whitebox{{\hbox{\hskip 1pt
 \vrule height 6pt depth 1.5pt
 \lower 1.5pt\vbox to 7.5pt{\hrule width
    3.2pt\vfill\hrule width 3.2pt}%
 \vrule height 6pt depth 1.5pt
 \hskip 1pt } }}
\def\qed{\ifhmode\allowbreak\else\nobreak\fi\hfill\quad\nobreak
     \whitebox\medbreak}

\newcommand{\ignore}[1]{}

\theoremstyle{plain}
\newtheorem{theorem}{Theorem}[section]

\newtheorem{lemma}[theorem]{Lemma}
\newtheorem{definition}[theorem]{Definition}
\newtheorem{problem}[theorem]{Problem}
\newtheorem{proposition}[theorem]{Proposition}
\newtheorem{example}[theorem]{Example}

\newtheorem{remark}[theorem]{Remark}

\def\qed{{\hfill$\square$}}
\def\proof{{\vspace{-0.3cm}\bf Proof: \,}}

\def\Z{{\Bbb Z}}

\def\R{{\Bbb R}}
\def\T{{\Bbb T}}
\def\P{{\Bbb P}}
\def\F{{\Bbb F}}
\def\U{{\Bbb U}}
\def\mod{{\mathrm{mod\,\,}}}

\title{Divisible difference families 
from Galois rings $GR(4,n)$ and Hadamard matrices}
\author{Koji Momihara\footnotemark[1] \quad Mieko Yamada\footnotemark[2]}
\date{} 
\begin{document}
\maketitle
\footnotetext[1]{
Department of Mathematics, Faculty of Education, Kumamoto University,  
2-40-1 Kurokami, Kumamoto 860-8555, Japan; Email: 
momihara@educ.kumamoto-u.ac.jp}
\footnotetext[2]{
Institute of Science and Engineering, Kanazawa University, Kakuma-machi, Kanazawa, 920-1192, Japan; 
Email: myamada@se.kanazawa-u.ac.jp}
\renewcommand{\thefootnote}{\arabic{footnote}}
%%%%%%%%%%%%%%%%%%%%%%%%%%%%%%%%%%%%%%%%%%%%%%%%%%%%%
%%%%%%%%%%%%%%%%%%%%%%%%%%%%%%%%%%%%%%%%%%%%%%%%%%%%%
\begin{abstract}
We give a new construction  of difference families generalizing Szekeres's difference families~\cite{Sze}. 
As an immediate  consequence,  we obtain some new examples of  difference families with several blocks in  multiplicative subgroups of finite fields.
We also prove that there exists an infinite family of divisible difference families  with two blocks 
in a unit subgroup of the Galois ring $GR(4,n)$. Furthermore, we obtain a new construction method of 
symmetric Hadamard matrices by using 
divisible difference families and a new array. 
\end{abstract}
\begin{center} 
{\small Keywords: difference family; divisible difference family; Hadamard matrix}
\end{center}
%%%%%%%%%%%%%%%%%%%%%%%%%%%%%%%%%%%%%%%%%%%%%%%%%%%%%
%%%%%%%%%%%%%%%%%%%%%%%%%%%%%%%%%%%%%%%%%%%%%%%%%%%%%
\section{Introduction}
Let $(G,\cdot)$ be a finite abelian group with identity $1$ and let $N$ be 
a subgroup of $G$. 
Let $D_i$ be a $k_i$-subset (called a {\it block}) of $G$ for $i = 1,2,...,b$,  and set ${\mathcal F}=\{D_i\,|\,1\le i\le b\}$. 
For a block $D_i$, let $\theta_i(d) = |\{(x, y)|x \cdot y^{-1}, x \ne y , x, y \in D_i \}|, i = 1,2,\ldots,b$.
Furthermore, we put $\theta_{\mathcal F}(d) = \sum_{i=1}^{b} \theta_i (d)$. 
%
%%%%%%%%%%%%%%%%%%%%%%%%%%% definition %%%%%%%%%%%%%%%%%%%%%%%%%%%%%%%%%%
%
\begin{definition}
A family ${\mathcal F}=\{D_i\,|\,1\le i\le b \}$ is called a {\it $(G,N, \{k_1,...,k_b \},\lambda,\mu)$ 
divisible difference family}  if
\[ \theta_{{\mathcal F}}(d) = \left\{ 
\begin{array}{ll}
\lambda &\quad  \text{if $d \in N \setminus \{1\}$,} \\
\mu &\quad \text{if $d \in G \setminus N$.}
\end{array}
\right. 
\]
\end{definition}%%%%%%%%%%%%%%%%%%%%%%%%%%%%%%%%%%%%%%%%%%%%%%%%%
In what follows, we abbreviate a divisible difference family as DDF.
If $a_i$ blocks have same cardinality $k_i$, we write $k_i^{(a_i)}$ simply. If $k=k_1 = k_2 = \cdots = k_b$, 
then we will write $(G,N,k,\lambda, \mu)$-DDF. 
Especially if $N = \{1\}$, it is called a {\it $(G, \{k_1,\ldots,k_b \}, \lambda)$ difference family}, briefly denoted as DF.
If $b = 1$, then it is called a {\it $(G, N, k,\mu, \lambda)$ divisible difference set}. A divisible 
difference set with $N = \{1 \}$ is exactly a $(G, k, \lambda)$ difference set.

A divisible difference family yields a ``group divisible design'' by developing 
its blocks.  
Let $V$ be a set of $vg$ points, ${\mathcal B}$ be a collection of $k$-subsets ({\it blocks})
of $V$, and ${\mathcal G}$ be a collection of $g$-subsets (called {\it groups}) partitioning $V$. 
The triple $(V,{\mathcal B},{\mathcal G})$ is called a {\it $(k,\lambda,\mu)$-group divisible design of type $g^{v}$} if any two elements from different (resp. same) groups appear exactly $\mu$ (resp. $\lambda$) blocks. If 
$g=1$, then the pair $(V,{\mathcal B})$ is called a $2$-$(v,k,\lambda)$ design. 
For any block $B$
of a group divisible design $(V,{\mathcal B},{\mathcal G})$ and a permutation $\sigma$ on $V$,
define $B^\sigma=\{b^\sigma\,|\,b\in B\}$. If $B^\sigma\in {\mathcal B}$ for all
$B\in {\mathcal B}$, then $\sigma$ is called an {\it automorphism} of the design.
The set of all such permutations forms a group under composition, which is  called the
{\it full automorphism group}. Any of its subgroups is called an
{\it automorphism group}. 
It is obvious that a $(G,N,k,\lambda,\mu)$-DDF ${\mathcal F}$ yields a 
$(k,\lambda,\mu)$ group divisible design of type $|N|^{|G|/|N|}$ by developing its blocks, i.e., ${\mathcal B}=\{D\cdot x\,|\,D\in {\mathcal F},x\in G\}$ with $G$ as its automorphism group. Divisible difference families (sets) 
and difference families (sets)
have been extensively studied as a major topic in the design theory in relation to group divisible designs and $2$-designs. We refer the reader to \cite{AHJP,Arasu,Co,Davis,H,Ionin,Kan,J1,J2,Sarvate,WH} on divisible difference sets and group divisible designs. However, there are not so many papers about divisible difference families with a few blocks. 

%The divisible difference sets have been studied and we see the recent remarkable work in
%\cite{Fer,Ionin,Kan,WH}. On the other hand, the study of the constuction of the divisible
%difference family with at least two blocks and large blocksize has not been investigated.
%We gives the references \cite{Arasu,Sarvate}.

In \cite{Sze},  Szekeres constructed difference families with two blocks in finite fields as follows.
%
%%%%%%%%%%%% propsition %%%%%%%%%%%%%%%%%%
%
\begin{proposition}\label{SDSSze}
Let $q\equiv 3\,(\mod{4})$ be a prime power.  
Let $\F_q$ be the finite field with $q$ elements and $N$ be the set of nonzero squares in 
$\F_q$. Then, the set of 
\[
D_1:=(N-1)\cap N \mbox{ and }  D_2:=(N+1)\cap N
\]
forms an $(N,(q-3)/4,(q-7)/4)$-DF. 
\end{proposition}%%%%%%%%%%%%%%%%%%%%%%%%%%%%%%%%%%%
We say that the difference family of Proposition~\ref{SDSSze} is the 
{\it Szekeres difference family}. The  
Szekeres difference family $\{D_1, D_2 \}$ satisfies that  
\begin{equation}\label{skewness} 
\mbox{$-a\not\in D_i$ if $a\in D_i$ for either $i=1$ or $2$.}
\end{equation}

If an $n \times n$-matrix $M$ with entries from $\{ 1, -1 \}$ satisfies $MM^T = nI_n$, $M$ is
called an {\it Hadamard matrix}, where $I_n$ is the identity matrix of order $n$. An Hadamard matrix
$M$ is {\it symmetric} or {\it skew} if 
$M=M^T$ or $M-I_n=-(M-I_n)^{T}$, respectively. For a general background on Hadamard matrices, we refer to \cite{Hedayat,Horadam}.  There have been known several constructions of Hadamard matrices from 
difference families \cite{Seberry, Sze2,XX,XXS}. (These authors used the terminology  ``supplementary difference sets'' instead of difference families.) 
For example, if a $(G, m, m-1)$-DF with $|G| = 2m+1$ has the property (\ref{skewness}), we are able to  construct a skew Hadamard matrix as follows.
%
%%%%%%%%%%%%%%%%%%%% proposition %%%%%%%%%%%%%%%%%%%%%%%%%%%%%%
%
\begin{proposition}\label{SzeHada}
{\rm (Theorem 4.4  of p.~321 \cite{Sze2})}
Let ${\mathcal F}=\{D_1,D_2\}$ be a $(G,m,m-1)$-DF with $|G|=2m+1$ satisfying the condition (\ref{skewness}). 
We define $(2m+1)\times (2m+1)$-matrices $A=(a_{i,j})$ and $B=(b_{i,j})$ by 
\begin{align*}
&a_{i,j}:=2f_{D_1}(i-j)-1\mbox{\, \, for the characteristic function $f_{D_1}$ of $D_1$, and }\\
&b_{i,j}:=2f_{D_2}(i+j)-1\mbox{\, \, for the characteristic function $f_{D_2}$ of $D_2$},  
\end{align*}
where rows and columns $i,j$ are labeled by the elements of $G$. 
Then,  
the matrix $M$ defined by 
\begin{eqnarray}\label{array1}
M=\left[
\begin{array}{cccc}
1&1&{\bf e}^T&{\bf e}^T\\
-1&1&{\bf e}^T&-{\bf e}^T\\
-{\bf e}&-{\bf e}&A&B \\
-{\bf e}&{\bf e}&-B&A
 \end{array}
\right]
\end{eqnarray}
forms a skew Hadamard matrix of order $4(m+1)$, where 
${\bf e}$ is a column vector of length $2m+1$ with all one entries.  
\end{proposition}%%%%%%%%%%%%%%%%%%%%%%%%%%%%%%%%%%%%%
It is well known   that the set $N$ of Proposition~\ref{SDSSze} is a 
difference set in the additive group of $\F_q$, which yields an Hadamard matrix of order $q+1$. One can see that this Hadamard matrix is equivalent to that obtained 
by applying Proposition~\ref{SzeHada} to the difference familiy of 
Proposition~\ref{SDSSze}. 

%Thus, we can recover Hadamard matrices from the Szekeres difference families. (It is well known   that the set $N$ of Proposition~\ref{SDSSze} is a 
%difference set in the additive group of $\F_q$, which yields an Hadamard matrix of order $q+1$.) 

The main objective of this paper is to obtain an analogy in Galois rings of the relation between difference families and Hadamard matrices of Propositions~\ref{SDSSze} and \ref{SzeHada}. This paper is organized as follows:
In Section~\ref{mainconst}, we give a general construction method of difference families by 
generalizing Szekeres's construction of Proposition~\ref{SDSSze}. 
Furthermore, we obtain some new examples of difference families with several blocks in multiplicative subgroups of 
finite fields. In Section~\ref{sec:Gal}, we prove that there exists an infinite family of divisible difference families with two blocks 
in unit groups of Galois rings of characteristic four. Finally, in Section~\ref{sec:Had}, we give a new construction of 
symmetric Hadamard matrices by using a divisible difference family satisfying certain conditions  and an array obtained by modifying the matrix $M$ of  (\ref{array1}). 
\section{Generalization of  Szekeres's construction}\label{mainconst}
At first of this section, we give a general construction of divisible difference families in a unit subgroup of a commutative ring $\R$ under the assumption of the existence of difference families in the additive group of $\R$.
%%%%%%%%%%% Lemma 2.1 %%%%%%%%%%%%%%%%%%%%%%%%%%%%%%%%%%%%%%
\begin{lemma}\label{them1}
Let $\R$ be a commutative ring with identity $1$. Let $\R^+$ be a additive group 
and $\R^\ast$ be a unit group of $\R$, and $I := \R \setminus \R^\ast$. 
Let $N\le \R^\ast$ and let $S$ be a complete system of representatives of $\R^\ast/N$. 
Let ${\mathcal F}=\{D_i\,|\,1\le i\le b\}$ be an $(\R^{+},\{k_i\,|\,1\le i\le b\},\lambda)$-DF.
Assume that each $D_i$ is fixed by $N$, i.e., $xD_i=D_i$ for all $x\in N$.  
We define subsets $D_{i,y}$ of $N$ as follows:
\[
D_{i,y}:=y^{-1}(D_i-1)\cap N
\]
for $1\le i\le b$ and $y\in S$. We put ${\mathcal F}'=\{D_{i,y}\,|\,1\le i\le b;y\in S\}$. Then we have
\[ \theta_{{\mathcal F'}}(t) = \lambda - \lambda_t, \]
where
\[
\lambda_t:=\sum_{i=1}^{b}|D_i\cap (D_i-t+1)\cap (I+1)|. 
\]
If $\lambda_t$ 
is constant for all $t\in N\setminus \{1\}$, say $\nu$, 
then ${\mathcal F}'$ is an $(N, \{k_{i,y} |\,1\le i\le b;y\in S\}, \lambda - \nu)$-DF.
If $\lambda - \lambda_t$ is a constant value $\mu$ for $t \in N\setminus L$ and $\eta$ for 
$t \in L\setminus \{1\}$, where $L$ is a subgroup of $N$, then 
${\mathcal F}'$ is an $(N, L, \{k_{i,y} |\,1\le i\le b;y\in S\}, \eta, \mu)$-DDF.
\end{lemma}%%%%%%%%%%%%%%%%%%%%%%%%%%%%%%%%%%%%%%%%%%%%%%%%%%%%%%%%%%%%%
\proof  
The multiplicity $\theta_{{\mathcal F'}}(t)$ for each $t\in N$ is given by 
\begin{eqnarray*}
\sum_{i=1}^{b}\sum_{y\in S}|D_{i,y}\cap t D_{i,y}|
&=&\sum_{i=1}^{b}\sum_{y\in S}|(D_i-1)\cap t (D_i-1)\cap yN|\\
&=&\sum_{i=1}^{b}|D_i\cap (D_i-t+1)\cap (\R^\ast+1)|\\
&=&\sum_{i=1}^{b}|D_i\cap (D_i-t+1)|-\sum_{i=1}^{b}|D_i\cap (D_i-t+1)\cap (I+1)|. 
\end{eqnarray*}
By the assumption that ${\mathcal F}$ is an $(\R^+,\{k_i\,|\,1\le i\le b\},\lambda)$-DF, it holds
\[
\sum_{i=1}^{b}| D_i\cap (D_i-t+1)|=\lambda,
\]
and hence we obtain the assertion. \qed

\vspace{0.2cm}
We denote the multiplicative group and the additive group of the finite field $\F_q$ by 
$\F_q^\ast$ and $\F_q^{+}$, respectively, and a primitive element of $\F_q$ by $g$. Furthermore,  we denote 
the residue ring $\Z/s\Z$ of rational integers by $\Z_s$. 

Now, we provide several new examples of difference families in a multiplicative subgroup of  the finite field $\F_q$ by applying 
Lemma~\ref{them1} to known ``cyclotomic'' difference sets. 
We say that a difference set $D$ in $\F_q^+$ is {\it cyclotomic} if 
$D$ is a subgroup of  $\F_q^\ast$ or the  union of a subgroup of  $\F_q^\ast$ and $\{0\}$. 
%
%%%%%%%%%%%%%%%%%%% proposition 2.2 %%%%%%%%%%%%%%%%%%%%%%%%%
%  
\begin{proposition}\label{exa1}
\begin{itemize}
\item[(i)] For any prime power  $q\equiv 3\,(\mod{4})$, there exists a 
$(\Z_{(q-1)/2},(q-3)/4,(q-7)/4)$-DF. 
\item[(ii)] For any prime power $q=1+4t^2$ with $t\equiv 1\,(\mod{2})$, there exists a 
$(\Z_{(q-1)/4},(q-5)/16,(q-21)/16)$-DF. 
\item[(iii)] For any prime power $q=9+64a^2=1+8b^2$ with $a\equiv b\equiv 1\,(\mod{2})$,  
there exists a 
$(\Z_{(q-1)/8},(q-9)/64,(q-73)/64)$-DF. 
\end{itemize} 
\end{proposition}%%%%%%%%%%%%%%%%%%%%%%%%%%%%%%
\proof
It is known that the multiplicative subgroup $D$ of index $e$ of $\F_q^\ast$ forms a 
cyclotomic $(q,|D|,\lambda)$ difference set in the following cases: 
(i) $e=2$ and $q\equiv 3\,(\mod{4})$; (ii)  $e=4$ and $q=1+4t^2$ with $t\equiv 1\,(\mod{2})$; and (iii) $e=8$ and $q=9+64a^2=1+8b^2$ with $a\equiv b\equiv 1\,(\mod{2})$. (See \cite{Storer}.) 
Now, we apply Lemma~\ref{them1} to these three difference sets $D$ as 
$N=D$. 
It is clear that $D$ is fixed by $N$.  
Let $S = \{g^i\,|\,0\le i\le e-1\}$ and $D_{1,g^i} = g^{-i}(D - 1) \cap D$ for $0\le i\le e-1$. 
For $t \in D \setminus \{1 \}$, we  have 
\[ 
\lambda_t = |D \cap (D - t +1) \cap (\{0\} + 1)| = 1. 
\]
Hence, by Lemma~\ref{them1}, $\{ D_{1,g^i}\,|\, 0\le i\le e-1\}$ forms a $(D,\{k_{1,y}\,|\,y\in S\},\lambda-1)$-DF. Since $N\simeq \Z_{(q-1)/e}$, this DF is isomorphic to that with the parameters in the statement of the proposition.  
Finally, we compute each block size $k_{1,y}=|D_{1,y}|$. 
By noting that $|(yN+1)\cap N|=|(N+y)\cap N|$ (cf. \cite[Lemma 3 (b)]{Storer}), 
we have 
\[
|D_{1,y}| =|y^{-i}(D - 1) \cap D|=|(yN+1)\cap N|=|(N+y)\cap N|=\lambda 
\]
since $D=N$ is a difference set in $\F_q^+$. 
This completes the proof. 
\qed

\vspace{0.2cm}
%
%%%%%%%%%%%%%%%%%%%%%%%%% proposition 2.3 %%%%%%%%%%%%%%%%%%%%%%%%%%
%%
\begin{proposition}\label{exa2}
\begin{itemize}
\item[(i)] For any prime power $q\equiv 3\,(\mod{4})$, there exists a $(\Z_{(q-1)/2},K,(q-3)/4)$-DF, 
where $K=\{(q-3)/4,(q+1)/4\}$. 
\item[(ii)] For any prime power $q=9+4t^2$ with $t\equiv 1\,(\mod{2})$, there exists a 
$(\Z_{(q-1)/4},K,(q-13)/16)$-DF, where $K=\{(q-13)/16,(q+3)/16^{(3)}\}$. 
\item[(iii)] For any prime power  $q=441+64a^2=49+8b^2$ with $a\equiv b\equiv 1\,(\mod{2})$, there exists a 
$(\Z_{(q-1)/8},K,(q-56)/64)$-DF, where $K=\{(q-56)/64,(q+8)/64^{(7)}\}$. 
\end{itemize} 
\end{proposition}%%%%%%%%%%%%%%%%%%%%%%%%%%%%%%%%%%%%%%%%%%%%%%%%%%%%%%%%%%%%
\proof
It is known that for a multiplicative subgroup $N$ of index $e$ of $\F_q^\ast$,  $D=N\cup \{0\}$ forms a 
cyclotomic $(q,|D|,\lambda)$ difference set in the following cases: 
(i) $e=2$ and $q\equiv 3\,(\mod{4})$; (ii)  $e=4$ and $q=9+4t^2$ with $t\equiv 1\,(\mod{2})$; and (iii) $e=8$ and $q=441+64a^2=49+8b^2$ with $a\equiv b\equiv 1\,(\mod{2})$. (See \cite{Storer}.) 
Now, we apply Lemma~\ref{them1} to these three difference sets $D$. 
It is clear that $D$ is fixed by $N$.  
Let $S = \{g^i\,|\,0\le i\le e-1\}$ and $D_{1,g^i} = g^{-i}(D - 1) \cap N$ for $0\le i\le e-1$. 
For $t \in D \setminus \{1 \}$, we  have 
\[ 
\lambda_t = |D \cap (D - t +1) \cap (\{0\} + 1)| = 1. 
\]
Hence, by Lemma~\ref{them1}, $\{D_{1,g^i}\,|\, 0\le i\le e-1\}$ forms a $(D,\{k_{1,y}\,|\,y\in S\},\lambda-1)$-DF. Since $N\simeq \Z_{(q-1)/e}$, this DF is isomorphic to that with the parameters in the statement of the proposition.  
Finally, we compute each block size $k_{1,y}=|D_{1,y}|$. 
By noting that $|(yN+1)\cap N|=|(N+y)\cap N|$, 
we have 
\begin{eqnarray*}
|D_{1,y}|&=&|y^{-1}(D - 1) \cap N|=|(yN + 1) \cap D|\\
&=&|(N+y)\cap N|+|(yN+1)\cap \{0\}|\\
&=&|(D+y)\cap D|-|\{y\}\cap D|-|(N+y)\cap \{0\}|+|\{-1\}\cap yN|. 
\end{eqnarray*}
Since $D$ is a difference set in $\F_q^+$, we have 
 $|(D+y)\cap D|=\lambda$. Furthermore, we have  
$-|\{y\}\cap D|-|(N+y)\cap \{0\}|+|\{-1\}\cap yN|=-1$ or $0$ according to $y\in N$ or not. This completes the proof. \qed

\vspace{0.2cm}
\begin{remark}
\begin{itemize}
\item[(1)] Let $q\equiv 3\,(\mod{4})$ and $N$ be the multiplicative subgroup of index $2$ of $\F_q^\ast$. Proposition~\ref{SDSSze} says that 
the set of $(N-1)\cap N$ and  $(N+1)\cap N$
forms an $(N,(q-3)/4,(q-7)/4)$-DF. On the other hand, 
Proposition~\ref{exa1}~(i) says that  
the set of $(N-1)\cap N$ and  $-(N-1)\cap N$
forms an $(N,(q-3)/4,(q-7)/4)$-DF. 
By noting that   
\begin{eqnarray*}
x\in ((N+1)\cap N)^{-1}&\Longleftrightarrow &
x\in N^{-1} \mbox{\, and\, }x\in (N+1)^{-1} \\
&\Longleftrightarrow &x\in N \mbox{\, and \, }x^{-1}-1\in N\\
&\Longleftrightarrow& x\in N \mbox{\, and \, }1-x=x(x^{-1}-1)\in N\\
&\Longleftrightarrow& x\in -(N-1)\cap N, 
\end{eqnarray*}
i.e., $ ((N+1)\cap N)^{-1}= -(N-1)\cap N$, 
we can claim that the Szekeres difference families of Proposition~\ref{SDSSze} and 
our difference families in Proposition~\ref{exa1} (i) are equivalent. 
\item[(2)] Let ${\mathcal F}$ be the $(\Z_{\frac{q-1}{e}},\{\lambda,(\lambda+1)^{(e-1)}\},\lambda)$-DF obtained in 
Proposition~\ref{exa2}. Consider the development of ${\mathcal F}$, i.e., 
${\mathcal B}=\{D+x\,|\,D\in {\mathcal F},x\in \Z_{\frac{q-1}{e}}\}$. 
By adding a new point $\infty$ to each block of size $\lambda$ in ${\mathcal B}$, we obtain a ($1$-rotational) $2$-$((q-1)/e+1,\lambda+1,\lambda)$ design on $\Z_{\frac{q-1}{e}}\cup \{\infty\}$. 
\end{itemize}
\end{remark}
%%%%%%%%%%%%%%%%%%%%%%%%%%%%%%%%%%%%%%%%%%%%%%%%%%%%%%%%%%%%%%%%%%%%%%%%%%%%%
%%%%%%%%%%%%%%%%%%%%%%% section 3 %%%%%%%%%%%%%%%%%%%%%%%%%%%%%%%%%%%%%
%%%%%%%%%%%%%%%%%%%%%%%%%%%%%%%%%%%%%%%%%%%%%%%%%%%%%%%%%%%%%%%%%%%%%%%%%%
\section{Divisible difference families from Hadamard difference sets in Galois rings $GR(4,n)$}\label{sec:Gal}
Let $q=p^r$ be a power of a prime $p$.  We denote the absolute trace from $\F_q$ to $\F_p$
by ${\rm Tr}_{\F_q}$. An {\it additive character} of $\F_q$ is written as
\begin{equation}\label{chara'}
\chi_b(c)= \zeta_p^{{\rm Tr}_{\F_{q}}(bc)} \mbox{ for any $b\in \F_q$},  
\end{equation}
where $\zeta_p={\rm exp}(\frac {2\pi i}{p})$.
Let $g(x) \in \Z_{p^2} [x]$ be a primitive basic irreducible polynomial of degree $n$ 
and denote the root of $g(x)$ by $\xi$. 
Then $\Z_{p^2} [x]/g(x)$ is  called 
a {\it Galois ring} of characteristic $p^2$ and of an extension degree $n$, and denoted by $GR(p^2, n)$. 
The algebraic extension of $\Z_{p^2}$ obtained by adjoining $\xi$ is isomorphic to 
$\Z_{p^2} [x]/g(x)$.
For easy reference, we put $\R_n = GR(4, n)$ as $p=2$. 
$\R_n$ has a unique maximal ideal  $\P_n =  p\R_n$ and the residue ring $\R_n/\P_n$ is 
isomorphic to $\F_{p^n}$. We take $\T_n=\{0,1,\xi,\ldots,\xi^{p^n-2}\}$ as a set of representatives of
 $\R_n/\P_n$. An arbitrary element $\alpha\in \R_n$ is uniquely represented as $\alpha = \alpha_0+p\alpha_1$, 
 $\alpha_0 ,\alpha_1\in \T_n$.

We define the map $\tau:\R_n \longrightarrow \T_n^\ast(:=\T_n\setminus \{0\})$ as $\tau(\alpha) = \alpha^{p^n}$.
The kernel of $\tau$ is the group $\U_n$ of principal units, which are elements of the form $1 + 2\beta,
\beta \in \T_n$. For the element $1 + 2\beta \in \U_n$, we may regard $\beta$ as an element of $\F_{p^n}$, 
then $\U_n$ is isomorphic to the additive group of of $\F_{p^n}$. 
If we denote the set of all units in $\R_n$ by $\R_n^\ast$, then $\R_n^\ast$ is the direct product of 
$\U_n$ and the cyclic group $\T_n^\ast$ of order $p^n - 1$. In other words, every element of 
$\R_n^\ast$ is uniquely 
represented as $\alpha_0(1+p\alpha_1)$, $\alpha_0, \alpha_1\in \T_n, \alpha_0 \not= 0$.  

In what follows, we identify $\R_n/p\R_n$ with $\F_q$ and denote each element of $\R_n/p\R_n$ and $\F_q$ by $\overline{a}$ 
in common. 
In this section, we construct divisible difference families in Galois rings $\R_n=GR(4,n)$  by  applying Lemma~\ref{them1} to 
 difference sets obtained in \cite{YY}.
Let 
$E:=\{\overline{x}\in \F_{2^n}\,|\,\mathrm{Tr}_{\F_{2^{n}}}(\overline{ux})=\bar{0}\}$ for 
$\overline{u}\in \F_{2^n}^\ast$ such that $\mathrm{Tr}_{\F_{2^{n}}}(\overline{u})=\overline{0}$. Note that 
$E$ is a subgroup of order $2^{n-1}$ of $\F_{2^n}^+$. 
In \cite{YY}, it was proved that $D=\{a(1+2b)\,|\,a\in \T_n^\ast,b\in \T_n\mbox{ such that }
\overline{b}\in E)\}$ forms an $(\R_n^{+},2^{n-1}(2^n-1),2^{n-1}(2^{n-1}-1))$ difference set. (A difference set with parameters $(G,2^{n-1}(2^n-1),2^{n-1}(2^{n-1}-1))$ for a group $G$ of order $2^{2n}$ is called {\it Hadamard}.)
Notice that $D$ is a subgroup of index $2$ of $\R_n^\ast$ isomorphic to $\Z_{2^n-1}\times \Z_{2}^{n-1}$.
We shall use the fact that the characteristic function $\psi  (\overline{x})$ of $E$ 
in $\F_{2^{n}}$ is given by  
\[
\psi (\overline{x})=\frac{1}{2}\sum_{\overline{w}\in \F_2}\chi(\overline{uwx})
=\frac{1}{2}\sum_{\overline{w}\in \F_2}(-1)^{\mathrm{Tr}_{\F_{2^{n}}}(\overline{uwx})},  
\]
where $\chi$ is the canonical additive character of $\F_{2^n}$. 
%%%
%%%%%%%%%%%%%%%%%%%% theorem 3.1 %%%%%%%%%%%%%%%%%%%%%%%%%%
%%%%
\begin{theorem}\label{mainthem}
Let $D$ and $E$ be as in the above. Let $N$ be a subgroup of $D$ and put $L=N\cap \U_n$.
Then there exists an $(N,L,K,2^{n}(2^{n-2}-1),2^{n-1}(2^{n-1}-1)-2^{n-2})$-DDF
where $K$ is the set of cardinalities of blocks. In particular, if $N=D$, 
the DDF is isomorphic to a   
$(\Z_{2^n-1}\times \Z_{2}^{n-1},\{0\}\times \Z_{2}^{n-1},2^{n-1}(2^{n-1}-1),2^{n}(2^{n-2}-1),
2^{n-1}(2^{n-1}-1)-2^{n-2})$-DDF. 
\end{theorem}%%%%%%%%%%%%%%%%%%%%%%%%%%%%%%%%%%%%%%%%%%%%%%%%%%%%%%%%%%%%%%%%%%%%%
\proof
Obviously, $D$ is fixed by $N$ since $N\le D$. Furthermore, we have $D \cap (\P_n + 1) = \{1 +2a |\,\overline{a}\in E\}$. 
Let $t=c(1+2d)\in N$. 
By applying Lemma~\ref{them1} to $D$, the number 
$\lambda_t$ is computed  as follows:
\begin{eqnarray*}
\lambda_t&=&|D\cap (D-t+1)\cap (\P_n+1)|\\
&=&|\{1+2a\,|\,\overline{a}\in E\}\cap (\{c'(1+2a')\,|\,c'\in \T_n^\ast,\overline{a'}\in  E\}-c(1+2d)+1)|\\
&=&|E\cap \{\overline{c(a'-d)}\,|\,\overline{a'}\in  E\}|\\
&=&\sum_{\overline{x}\in \F_{2^n}}\Big(\frac{1}{2}\sum_{\overline{w}\in \F_2}\chi(\overline{uwx})\Big)\cdot 
\Big(\frac{1}{2}\sum_{\overline{v}\in \F_2}
\chi(\overline{uv}(\overline{c}^{-1}\overline{x}+\overline{d}))\Big)\\
&=&\frac{1}{4}\sum_{\overline{x}\in \F_{2^n}}\Big(\sum_{\overline{w}\in \F_2}
\sum_{\overline{v}\in \F_2}\chi(\overline{ux}(\overline{w}+\overline{c}^{-1}\overline{v})+\overline{udv})\Big)\\
&=&
\left\{
\begin{array}{ll}
2^{n-2}+2^{n-2}\chi(\overline{ud})&  \mbox{if $\overline{c}^{-1}=\overline{1}$,}\\
2^{n-2}&  \mbox{otherwise,}
 \end{array}
\right. \hspace{5.5cm}
\\
&=&\left\{
\begin{array}{ll}
2^{n-1}&  \mbox{if $c=1$, that is, $t \in L$,}\\
2^{n-2}&  \mbox{otherwise. }
 \end{array}
\right.
\end{eqnarray*}
The set $K$ of cardinalities of blocks is determined by the subgroup $N$.
In particular, if $N=D$, the size $k_y$ of each block $y^{-1}(D-1)\cap N$ is given by 
\begin{eqnarray*}
k_y=|(yD+1)\cap D|
=|(D+y)\cap D|=2^{n-1}(2^{n-1}-1)
\end{eqnarray*}
since $|(yD+1)\cap D|=|(D+y)\cap D|$ and $D$ is a difference set in $\R_n^+$. 
Since $D\le \R_n^\ast$ is isomorphic to $ \Z_{2^n-1}\times \Z_{2}^n$, this DDF is 
isomorphic to that with the parameters in  the statement of  the theorem.
\qed
%
%%%%%%%%%%%%%%%%%%%%%%% Remark 3.2 %%%%%%%%%%%%%%%%%%%%%%%%%%%%%%%%%%%%%%
%

\vspace{0.2cm}
\begin{remark}
We notice that $D\cup \P_n$ forms an $(\R_n^{+},2^{n-1}(2^{n}+1),2^{n-1}(2^{n-1}+1))$ difference set. 
Let $N$ and $L$ be as in Theorem~\ref{mainthem}. We have 
an $(N,L,K,2^{2(n-1)},2^{n-2}(2^{n}+1))$-DDF 
since 
\[
|(D\cup \P_n)\cap ((D\cup \P_n)-t+1)\cap (\P_n+1)|=|D\cap (D-t+1)\cap (\P_n+1)|
\]
for any $t\in N\setminus \{1\}$. 
In particular, if $N=D$, 
it is isomorphic to a DDF with parameters 
$(\Z_{2^n-1}\times \Z_2^{n-1},\{0\}\times \Z_2^{n-1},2^{2(n-1)},2^{2(n-1)},2^{n-2}(2^{n}+1))$ 
since 
\begin{eqnarray*}
k_y&=&|y^{-1}(N\cup \P_n-1)\cap N|=|y^{-1}(N-1)\cap N|+|y^{-1}(\P_n-1)\cap N|\\
&=&|(N+y)\cap N|+|\P_n \cap (yN+1)|\\
&=&2^{n-1}(2^{n-1}-1)+2^{n-1}=2^{2(n-1)}. 
\end{eqnarray*}
\end{remark}
%
%%%%%%%%%%%%%%%%%%%%%%%% example 3.3 %%%%%%%%%%%%%%%%%%%%%%%%%%%%%%%%%%%%%%
%%%%%
\begin{example}
Let $\xi$ be a root of $g(x):=x^3+3x^2+2x+3\in \Z_4[x]$, a primitive basic irreducible polynomial of 
$GR(4,3)$. Furthermore, let $xyz$ denote the element $x\xi^2+y\xi+z\in GR(4,3)$. 
Let $D=\{a(1+2b)\,|\,a\in \T_3^\ast ,b\in \{0,1, \xi^2,
\xi^3 \}\}$ and $L=D\cap \U_3$. Note that $D$ forms an $(\R_n^{+},28,12)$ difference set. Also, $D$ is a subgroup
of  index $2$ of $\R_3^\ast$.
Then, the family of 
\begin{eqnarray*}
D_1&=&\{103,232,322,112,211,111,231,121,300,332,212,331\} \mbox{ and}\\ 
D_2&=&\{233,322,332,113,213,121,010,333,103,300,112,030\} 
\end{eqnarray*}
forms a $(D,L,12,8,10)$-DDF.
Since  $\R_3^\ast \cong \Z_7\times \Z_2^3$ and $L \cong \{0\}\times \Z_2\times \Z_2$,
the blocks can be expressed as    
\begin{eqnarray*}
\mathcal{D}_1&=&\{(1,0,1),(1,1,1),(2,0,1),(2,1,0),(3,0,0),(3,0,1),\\
& &\, \, (4,0,0),(4,0,1),(5,0,1),(5,1,0),(6,0,1),(6,1,1)\}\subseteq \Z_7\times \Z_2^2 \mbox{ and}\\
\mathcal{D}_2&=&\{(1,0,0),(1,1,0),(2,0,1),(2,1,0),(3,0,0),(3,0,1),\\
& &\, \, (4,1,0),(4,1,1),(5,0,0),(5,1,1),(6,0,1),(6,1,1)\}\subseteq \Z_7\times \Z_2^2, 
\end{eqnarray*}
and the family $\{\mathcal{D}_1, \mathcal{D}_2 \}$ forms a $(\Z_7\times \Z_2\times \Z_2,\{0\}\times \Z_2\times \Z_2,12,8,10)$-DDF. 
\end{example}%%%%%%%%%%%%%%%%%%%%%%%%%%%%%%%%%%%%%%%%%%%%%%%%%%%%%%%%%%%%%%%%%%%
Next we show some important and interesting properties of divisible difference families
obtained in Theorem~\ref{mainthem}. 
%%%%%%%
%%%%%%%%%%%%%%%%%%%%% proposition 3.4 %%%%%%%%%%%%%%%%%%%%%%%%%%%%5
%%%%%%%%%%%
\begin{proposition}\label{interest}
The $(\Z_{2^n-1}\times \Z_{2}^{n-1},\{0\}\times \Z_{2}^{n-1},2^{n-1}(2^{n-1}-1),2^{n}(2^{n-2}-1),2^{n-1}(2^{n-1}-1)-2^{n-2})$-DDF 
${\mathcal F}=\{\mathcal{D}_1,\mathcal{D}_2\}$
of Theorem~\ref{mainthem} has the following properties: 
\begin{itemize}
\item[(i)] If $a\in {\mathcal D}_1$, then $-a\in {\mathcal D}_1$. 
\item[(ii)] If $a\in {\mathcal D}_2$, then $-a\not \in \mathcal{D}_2$. 
\item[(iii)] For each $i=1,2$, it holds that $|{\mathcal D}_i\cap (\{0\}\times \Z_{2}^{n-1})|=0$ and 
$|{\mathcal D}_i\cap (\{j\}\times \Z_{2}^{n-1})|=2^{n-2}$ for any $j\in\Z_{2^n-1}\setminus\{0\}$.  
\end{itemize}
\end{proposition}%%%%%%%%%%%%%%%%%%%%%%%%%%%%%%%%%%%%%%%%%%%%%%%
\proof
We consider the $(D,L,K,2^{n}(2^{n-2}-1),2^{n-1}(2^{n-1}-1)-2^{n-2})$-DDF ${\mathcal F}=\{D_1,D_2\}$ in $GR(4, n)$ obtained in Theorem~\ref{mainthem}, 
where $D_1=(D-1)\cap D$ and $D_2=y^{-1}(D-1)\cap D$. Here, $\{1,y\}$ is a set of representatives for $\R_n/D$. We can assume that $y\in \U_n$. Then, we have $y^2=1$.  

\vspace{0.2cm}
\noindent
{\bf (i)} We prove that  $a^{-1} \in (D-1)\cap D$  if $a \in (D-1)\cap D$.
\begin{eqnarray*}
& & a\in (D-1)\cap D\\
&\Longleftrightarrow& a\in D\mbox{ and }a+1\in D\\
&\Longleftrightarrow& a^{-1}\in D\mbox{ and }a^{-1}+1=a^{-1}(a+1)\in D\\
&\Longleftrightarrow& a^{-1}\in (D-1)\cap D. 
\end{eqnarray*}

\vspace{0.2cm}
\noindent
{\bf (ii)} We prove that $a^{-1}\not \in y(D-1)\cap D$ if $a\in y(D-1)\cap D$.
\begin{eqnarray*}
& & a\in y(D-1)\cap D\\
&\Longleftrightarrow& a\in D\mbox{ and }ay^{-1}+1\in D\\
&\Longleftrightarrow& a^{-1}\in D\mbox{ and }
a^{-1}y^{-1}+1=a^{-1}y+1=a^{-1}y(ay^{-1}+1)\in yD\\
&\Longleftrightarrow& a^{-1}\in y(yD-1)\cap D\\
&\Longleftrightarrow& a^{-1}\not \in y(D-1)\cap D. 
\end{eqnarray*}
Under the isomorphism $\phi:\R_n^\ast \longrightarrow \Z_{2^n-1}\times \Z_{2}^n$, we obtain 
the assertions (i) and (ii) in the theorem. 

\vspace{0.2cm}
\noindent
{\bf (iii)}  
The image of $\xi^j \{1+2x\,|\,\overline{x}\in E\}$ by the isomorphism $\phi$ is 
$\{j\}\times \Z_{2}^{n-1}$. Hence it is sufficient to show  
\begin{equation}\label{size}
(\gamma(j) :=)|y(D-1)\cap D\cap \xi^j \{1+2x\,|\,\overline{x}\in E\}|=
\begin{cases} 
0 & \quad \text{if $j=0,$} \\
2^{n-2} & \quad \text{if  $j \not= 0.$}
\end{cases}
\end{equation}
Assume that $j= 0$. Since $\{1+2x\,|\,\overline{x}\in E\}\subseteq D$, we have 
\[ 
\gamma(0) = |\{y(w(1+2x)-1)\,|\,w\in \T_n^\ast,\overline{x}\in E\}\cap \{1+2x\,|\,\overline{x}\in E\}|.
\]
Since $w \in \T_n^\ast$, $w-1 \notin \U_n$, so that $y(w-1) \notin \U_n$.  
For $w\in \T_n^\ast,\overline{x}\in E$,  we have 
$y(w(1+2x)-1) = y(w-1)+2ywx \notin \U_n$. It yields $\gamma(0) = 0$.   

Assume that $j \not= 0$. By letting $y = 1 + 2a$, we have
\[
\gamma(j) =|(w-1)+2w(x+a)-2a\,|\,w\in \T_n^\ast,\overline{x}\in E\}\cap \xi^j \{1+2x\,|\,\overline{x}\in E\}|.
\]
There is a unique $w\in \T_n^\ast\setminus\{1\}$ such that $w-1 \in \xi^j \U_n$.
We put $w-1 = \xi^j (1+2c)$ for such $w$. Then, we have  
\begin{eqnarray*}
\gamma(j) 
&=&|\{(\xi^j(1+2c)+1)(1+2x+2a)-1-2a\,|\,\overline{x}\in E\}\cap \xi^j \{1+2x\,|\,\overline{x}\in E\}|\\
&=&|\{2c\xi^j+(\xi^j+1)(2x+2a)-2a\,|\,\overline{x}\in E\}\cap \{2\xi^j x\,|\,\overline{x}\in E\}|\\
&=&|((\overline{1+\xi^{-j}})E+\overline{c}+\overline{a})\cap E|\\
&=&\sum_{\overline{x}\in \F_{2^n}}\Big(\frac{1}{2}\sum_{\overline{w}\in \F_2}\chi(\overline{uwx})\Big)\cdot 
\Big(\frac{1}{2}\sum_{\overline{v}\in \F_2}
\chi(\overline{uv}(\overline{x}-\overline{c}-\overline{a})(\overline{1+\xi^{-j}})^{-1})\Big)\\
&=&\frac{1}{4}\sum_{\overline{x}\in \F_{2^n}}\Big(\sum_{\overline{w}\in \F_2}
\sum_{\overline{v}\in \F_2}
\chi(\overline{ux}(\overline{w}+\overline{v}(\overline{1+\xi^{-j}})^{-1})+
\overline{u}\overline{v}(\overline{c}+\overline{a})(\overline{1+\xi^{-j}})^{-1})\Big)\\
&=&2^{n-2}.
\end{eqnarray*}
Here, we used that $\overline{1+\xi^{-j}}\not=\overline{0}$ and 
$\overline{1}+(\overline{1+\xi^{-j}})^{-1}\not=\overline{0}$ for any $j\not=0$.
Under the isomorphism $\phi:\R_n^\ast \longrightarrow \Z_{2^n-1}\times \Z_{2}^n$, we obtain 
the assertion (iii) of  the theorem.
\qed

\vspace{0.2cm}
The next proposition shows that  $(D+2)\cap \T_n^\ast$ forms a $(\T_n^\ast,2^{n-1}-1,2^{n-2}-1)$ difference set. Then it yields an Hadamard matrix
of order $2^n$.  
%%%%%%%%
%%%%%%%%%%%%%%%%%%%%%%%% proposition 3.5 %%%%%%%%%%%%%%%%%%%%%%%%%%%%%
%%%%%%%%%%%%
\begin{proposition}\label{subdes}
The set $(D+2)\cap \T_n^\ast$ forms a $(\T_n^\ast,2^{n-1}-1,2^{n-2}-1)$ difference set. 
\end{proposition}%%%%%%%%%%%%%%%%%%%%%%
\proof
We show that for any $\xi^\ell\in \T_n^\ast\setminus\{1\}$ 
\[
|(D+2)\cap \T_n^\ast\cap \xi^{\ell}((D+2)\cap \T_n^\ast)|=2^{n-2}-1. 
\]
Since $D=xD$ for $x\in \T_n^\ast$, we have 
\begin{eqnarray}
& &|(D+2)\cap \T_n^\ast\cap \xi^{\ell}((D+2)\cap \T_n^\ast)|\nonumber\\
&=&|(D+2)\cap (D+2\xi^\ell)\cap  \T_n^\ast|\nonumber\\
&=&|\{a(1+2b)+2\,|\,a,b\in \T_n,\overline{b}\in E\}\cap \{a(1+2b)+2\xi^\ell\,|\,a,b\in \T_n,\overline{b}\in E\}\cap  \T_n^\ast|. \label{main-comp} 
\end{eqnarray}
Here,  $a(1+2b)+2=a+2(ab+1)\in \T_n^\ast$ if and only if $a=b^{-1}$ since $a\in \T_n$. Similarly, $a(1+2b)+2\xi^\ell\in \T_n$ if and only if $a=b^{-1}\xi^\ell$. Thus, (\ref{main-comp}) is reformulated as
\begin{eqnarray*}  
&=& |\{b^{-1}(1+2b)+2\,|\,\overline{b}\in E\setminus\{\overline{0}\}, b\in \T_n\}\cap 
\{b^{-1}\xi^{\ell}(1+2b)+2\xi^\ell\,|\,\overline{b}\in E\setminus\{\overline{0}\}, b\in \T_n\}|\\
&=&|\{b^{-1}\,|\,\overline{b}\in E\setminus\{\overline{0}\}, b\in \T_n\}\cap 
\{b^{-1}\xi^\ell\,|\,\overline{b}\in E\setminus\{\overline{0}\}, b\in \T_n\}|\\
&=&|E\setminus\{\overline{0}\}\cap 
\overline{\xi^{-\ell}}(E\setminus\{\overline{0}\})|\\
&=&\sum_{\overline{x}\in \F_{2^n}^\ast}\Big(\frac{1}{2}\sum_{\overline{w}\in \F_2}\chi(\overline{uwx})\Big)\cdot 
\Big(\frac{1}{2}\sum_{\overline{v}\in \F_2}
\chi(\overline{uvx\xi^{\ell}})\Big)\\
&=&\frac{1}{4}\sum_{\overline{x}\in \F_{2^n}^\ast}\Big(\sum_{\overline{w}\in \F_2}
\sum_{\overline{v}\in \F_2}
\chi(\overline{ux}(\overline{w}+\overline{v}\overline{\xi^{\ell}}))\Big)\\
&=&2^{n-2}-1. 
\end{eqnarray*}
The size of $(D+2)\cap \T_n^\ast$ is computed  as follows:  
\begin{eqnarray*}
|(D+2)\cap \T_n^\ast|
&=&|\{a(1+2b)+2\,|\,a,b\in \T_n,\overline{b}\in E\}\cap \T_n^\ast|\\
&=&|\{b^{-1}(1+2b)+2\,|\,\overline{b}\in E\setminus \{\overline{0}\},b\in \T_n\}\cap \T_n^\ast|\\
&=&|\{b^{-1}\,|\,\overline{b}\in E\setminus\{\overline{0}\}, b\in \T_n\}\cap \T_n^\ast|\\
&=&|E|-1=2^{n-1}-1.
\end{eqnarray*}
This completes the proof. 
\qed

\vspace{0.2cm}
\section{Construction of symmetric Hadamard matrices}\label{sec:Had}
In this section, we give a new construction of symmetric Hadamard matrices. 
%%%%%%%%%
%%%%%%%%%%%%%%%%%%%%%%% theorem 4.1 %%%%%%%%%%%%%%%%%%%%%%%%%%%%%%%%%%%%%%%%%%%%%%%
%%%%%%
\begin{theorem}\label{Hadamard}
Assume that there exists an Hadamard matrix of order $n$. Let ${\mathcal F}=\{D_1,D_2\}$ be 
a $(G,N,\frac{n(n-2)}{4},\frac{n(n-4)}{4},\frac{n(n-3)}{4})$-DDF with 
$|G|=\frac{n(n-1)}{2}\mbox{ and }|N|=\frac{n}{2}$.

Further assume that ${\mathcal F}$ satisfies the the following conditions: 
\begin{itemize}
\item[(i)] If $a\in D_1$, then $-a\in D_1$; 
\item[(ii)] $|D_i\cap N|=0$ and 
$|D_i\cap (N+j)|=n/4$ for any $j\in G\setminus N$ and for each $i=1,2$.
\end{itemize}
Then, there exists a symmetric Hadamard matrix of order $n^2$. 
\end{theorem}%%%%%%%%%%%%%%%%%%%%%%%%%%%%%%%%%%%%%%%%%%%%%%%%
\proof
We define the matrices $A=(a_{i,j}), B=(b_{i,j})$ and $C=(c_{i,j})$ of size $|G|\times |G|$ by 
\begin{align*}
&a_{i,j}:=f_{D_1}(i-j)\mbox{\,\,  for the characteristic function $f_{D_1}$ of $D_1$},\\
&b_{i,j}:=f_{D_2}(i+j)\mbox{\,\,  for the characteristic function $f_{D_2}$ of $D_2$},\\
&c_{i,j}:=f_{N}(i+j)\mbox{\,\,  for the characteristic function $f_{N}$ of $N$},
\end{align*}
for $i,j\in G$. Further we let $A':=2A-J_{|G|}$ and $B':=2B-J_{|G|}$, where $J_{|G|}$ is 
the all one matrix of size $|G|\times |G|$. Without loss of generality, 
we may assume that the first column and first row of the assumed Hadamard matrix $H$ are all one vectors.
Let $H'$ be the $(n-1)\times n$-matrix obtained by removing the first row of  $H$. Assume that 
the rows of the matrix $H'$ are indexed by the elements of $S=G/N$ and set $H'=({\bf h}_{i})_{i\in S}$.
Let $H_1=({\bf h}_i)_{i\in G}$ be the $n(n-1)/2\times n$-matrix defined by 
${\bf h}_{i+x}:={\bf h}_i$ for $i\in S$ and $x\in N$. 
Furthermore, let $H_2=(-{\bf h}_{-i})_{i\in G}$, then $H_2$ has size $n(n-1)/2\times n$.
Then the matrix
\begin{eqnarray*}
M=\left[
\begin{array}{ccc}
-J_{n}&H_1^{T} &H_2^{T} \\
H_1&B'&A' \\
H_2&A'&-B'-2C 
 \end{array}
\right]
\end{eqnarray*}
forms the desired symmetric Hadamard matrix of order $n^2$. 
It is clear that 
{\tiny
\begin{align*}
&MM^T=
\begin{bmatrix}
-J_{n}&H_1^{T} &H_2^{T} \\
H_1&B'&A' \\
H_2&A'&-B'-2C 
\end{bmatrix}
\begin{bmatrix}
-J_{n}&H_1^{T} &H_2^{T} \\
H_1&B'^T &A'^T \\
H_2&A'^T &-B'^T-2C^T 
\end{bmatrix}\\
&=  
\begin{bmatrix}
nJ_n + H_1^T H_1 + H_2^T H_2 & -J_n H_1^T +  H_1^T B'^T + H_2^T {A'}_1^T 
& -J_n H_2^T + H_1^T {A'}^T - H_2^T B'^T - 2H_2^T C^T \\
-H_1 J_n + B' H_1 + A' H_2 & H_1 H_1^T + A'A'^T + B'B'^T & H_1 H_2^T + B'{A'}^T - A'{B'}^T - 2A'C^T \\
-H_2J_n + A'H_1 - B'H_2 -2CH_2 & H_2 H_1^T + A'B'^T - B'A'^T -2CA'^T 
& H_2 H_2^T + A'A'^T + B'B'^T + 4CC^T+2B'C^T +2C{B'}^T
\end{bmatrix}.
\end{align*}
}
To prove $MM^T = n^2 I$, we show the following claims. 

\vspace{0.2cm}
\noindent
%%%%%%%%%%%%% claim 1 %%%%%%%%%%%%%%%%%%%%%%% 
{\bf Claim 1.} Write $H_1=(h^{(1)}_{i,j})$ and $H_2=(h^{(2)}_{i,j})$, where $i\in G,1\le j\le n$. Then, for each $\ell=1,2$ 
\begin{eqnarray*}
(\mbox{the $(i,j)$ entry of }H_\ell H_\ell^T)
=\sum_{k=1}^nh^{(\ell)}_{i,k}h^{(\ell)}_{j,k}
=\left\{
\begin{array}{ll}
n&  \mbox{if $i-j\in N$,}\\
0&  \mbox{otherwise. }
 \end{array}
\right.
\end{eqnarray*}
%%%%%%%%%%%%%%%%% claim 2 %%%%%%%%%%%%%%%%%%%%%  

\vspace{0.2cm}
\noindent
{\bf Claim 2.}   For each $\ell=1,2$,  
\begin{eqnarray*}
(\mbox{the $(i,j)$ entry of }H_\ell^T H_\ell)
=\sum_{k\in G}h^{(\ell)}_{k,i}h^{(\ell)}_{k,j}
=\left\{
\begin{array}{ll}
n(n-1)/2&  \mbox{if $i=j$,}\\
-n/2&  \mbox{otherwise. }
 \end{array}
\right.
\end{eqnarray*}
%%%%%%%%%%%%%%%%%%%%Claim3%%%%%%%%%%%%%%%%%%%%%%%%%%%%%%%%%%%%

\vspace{0.2cm}
\noindent
{\bf Claim 3.} Note that $h^{(2)}_{i,j}=-h^{(1)}_{-i,j}$, by the definition of $H_2$. Then, 
\begin{eqnarray*}
(\mbox{the $(i,j)$ entry of }H_1 H_2^T)
&=&\sum_{k=1}^nh^{(1)}_{i,k}h^{(2)}_{j,k}\\
&=&-\sum_{k=1}^n h^{(1)}_{i,k}h^{(1)}_{-j,k}\\
&=&\left\{
\begin{array}{ll}
-n&  \mbox{if $i+j\in N$,}\\
0&  \mbox{otherwise. }
 \end{array}
\right.
\end{eqnarray*}
%%%%%%%%%%%%%%%%%%%%Claim4%%%%%%%%%%%%%%%%%%%%%%%%%%%%%%%%%%%%

\vspace{0.2cm}
\noindent
{\bf Claim 4.} 
\begin{eqnarray*} 
(\mbox{The $(i,j)$ entry of }C C^T)
&=&\sum_{k\in G}c_{i,k}c_{j,k}\\
&=&\sum_{x\in S}\sum_{y\in N}c_{i,x+y}c_{j,x+y}\\
&=&\sum_{x\in S}\sum_{y\in N}c_{i,x}c_{j,x}\\
&=&\frac{n}{2}\sum_{x\in S}c_{i,x}c_{j,x}\\
&=&\left\{
\begin{array}{ll}
n/2&  \mbox{if $i-j\in N$,}\\
0&  \mbox{otherwise. }
 \end{array}
\right.
\end{eqnarray*}
%%%%%%%%%%%%%%%%%% claim 5 %%%%%%%%%%%%%%%%%%%%%%%%%% 

\vspace{0.2cm}
\noindent
{\bf Claim 5.} 
 By the assumption that $A$ and $B$ are incidence matrices of $D_1$ and $D_2$, respectively, 
\begin{eqnarray*}
(\mbox{the $(i,j)$ entry of }AA^T+BB^T)=
\left\{
\begin{array}{ll}
\frac{n(n-2)}{2}&  \mbox{if $i=j$,}\\
\frac{n(n-4)}{4}&  \mbox{if $i-j\in N$ and $i\not=j$,}\\
\frac{n(n-3)}{4}&  \mbox{otherwise. }
 \end{array}
\right.
\end{eqnarray*}
Furthermore, 
by the definition of $A'$ and $B'$, 
\begin{eqnarray*}
(\mbox{the $(i,j)$ entry of }A'{A'}^T+B'{B'}^T)=
\left\{
\begin{array}{ll}
n(n-1)&  \mbox{if $i=j$,}\\
-n&  \mbox{if $i-j\in N$ and $i\not=j$,}\\
0&  \mbox{otherwise. }
 \end{array}
\right.
\end{eqnarray*}
%%%%%%%%%%%%%%%%%% claim 6 %%%%%%%%%%%%%%%%%%%%%%%%

\vspace{0.2cm}
\noindent
{\bf Claim 6.}
 Write $A'=({a'}_{i,j})$, where $i,j\in G$. 
Note that for $x\in S$ 
\[
\sum_{y\in N}{a'}_{i,x+y}=\left\{
\begin{array}{ll}
-n/2&  \mbox{if $x-i\in N$,}\\
0&  \mbox{if $x-i\not\in N$,}
 \end{array}
\right.
\]
by the assumption (ii).
\begin{eqnarray*}
(\mbox{The $(i,j)$ entry of }A'C^T)
&=&\sum_{k\in G}{a'}_{i,k}c_{j,k}\\
&=&\sum_{x\in S}c_{j,x}\sum_{y\in N}{a'}_{i,x+y}\\
&=&-\frac{nc_{j,i}}{2}=
\left\{
\begin{array}{ll}
-n/2&  \mbox{if $i+j\in N$,}\\
0&  \mbox{otherwise. }\\
 \end{array}
\right.
\end{eqnarray*}
%%%%%%%%%%%%%%%%%%% claim 7 %%%%%%%%%%%%%%%%%%%%%%%

\vspace{0.2cm}
\noindent
{\bf Claim 7.} 
Write $B'=({b'}_{i,j})$, where $i,j\in G$. 
Note that for $x\in S$ 
\[
\sum_{y\in N}{b'}_{i,x+y}=\left\{
\begin{array}{ll}
-n/2&  \mbox{if $x+i\in N$,}\\
0&  \mbox{if $x+i\not\in N$,}
 \end{array}
\right.
\]
by the assumption (ii). Then, 
\begin{eqnarray*}
(\mbox{the $(i,j)$ entry of }B'C^T)
&=&\sum_{k\in G}{b'}_{i,k}c_{j,k}\\
&=&\sum_{x\in S}c_{j,x}\sum_{y\in N}{b'}_{i,x+y}\\
&=&-\frac{nc_{j,-i}}{2}=
\left\{
\begin{array}{ll}
-n/2&  \mbox{if $i-j\in N$,}\\
0&  \mbox{otherwise. }\\
 \end{array}
\right.
\end{eqnarray*}
Similarly, we have 
\begin{eqnarray*}
(\mbox{the $(i,j)$ entry of }C{B'}^T)=
\left\{
\begin{array}{ll}
-n/2&  \mbox{if $i-j\in N$,}\\
0&  \mbox{otherwise. }\\
 \end{array}
\right.
\end{eqnarray*}
%%%%%%%%%%%%%%% claim 8 %%%%%%%%%%%%%%%%%%%%%%%%%

\vspace{0.2cm}
\noindent
{\bf Claim 8.}
\begin{eqnarray*}
(\mbox{The $(i,j)$ entry of }B'H_1+A'H_2)
&=&\sum_{k\in G}{b'}_{i,k}h^{(1)}_{k,j}+\sum_{k\in G}{a'}_{i,k}h^{(2)}_{k,j}\\
&=&\sum_{x\in S}h^{(1)}_{x,j}\sum_{y\in N}{b'}_{i,x+y}+\sum_{x\in S}h^{(2)}_{x,j}\sum_{y\in N}{a'}_{i,x+y}\\
&=&-nh^{(1)}_{-i,j}/2-nh^{(2)}_{i,j}/2=0. 
\end{eqnarray*}
%%%%%%%%%%%%%%%%%%% claim 9 %%%%%%%%%%%%%%%%%%%%%

\vspace{0.2cm}
\noindent
{\bf Claim 9.}
Note that for $x\in S$ 
\[
\sum_{y\in N}{c}_{i,x+y}=\left\{
\begin{array}{ll}
n/2&  \mbox{if $x+i\in N$,}\\
0&  \mbox{if $x+i\not\in N$.}
 \end{array}
\right.
\]
Then, 
\begin{eqnarray*}
& &(\mbox{the $(i,j)$ entry of }A'H_1-B'H_2-2CH_2)\\
&=&\sum_{k\in G}{a'}_{i,k}h^{(1)}_{k,j}-\sum_{k\in G}{b'}_{i,k}h^{(2)}_{k,j}-2\sum_{k\in G}{c}_{i,k}h^{(2)}_{k,j}\\
&=&\sum_{x\in S}h^{(1)}_{x,j}\sum_{y\in N}{a'}_{i,x+y}-\sum_{x\in S}h^{(2)}_{x,j}\sum_{y\in N}{b'}_{i,x+y}-
2\sum_{x\in S}h^{(2)}_{x,j}\sum_{y\in N}{c}_{i,x+y}\\
&=&-nh^{(1)}_{i,j}/2+nh^{(2)}_{-i,j}/2-nh^{(2)}_{-i,j}=0. 
\end{eqnarray*}
%%%%%%%%%%%%%%%% claim 10 %%%%%%%%%%%%%%%%%%% 

\vspace{0.2cm}
\noindent
{\bf Claim 10.} It holds that
\[
A'{B'}^{T}=B'{A'}^{T}
\]
since 
\begin{eqnarray*}
(\mbox{the $(i,j)$ entry of }A'{B'}^T)
&=&\sum_{k\in G}(2f_{D_1}(i-k)-1)(2f_{D_2}(j+k)-1)\\
&=&\sum_{h\in G}(2f_{D_1}(j-h)-1)(2f_{D_2}(i+h)-1)\\
&=&(\mbox{the $(i,j)$ entry of }B'{A'}^T),
\end{eqnarray*}
where we put $h=j-i+k$. 

\vspace{0.2cm}

Now, we prove $MM^T=n^2I$ using the claims above.  
From Claim 2, we have
\begin{equation*}
(\mbox{the $(i,j)$ entry of }nJ_n + H_1^T H_1 + H_2^T H_2) =
n + 
\begin{cases}
n^2 - n & \text{if $i=j$,}\\
-n & \text{otherwise,}
\end{cases}
=
\begin{cases}
n^2 & \text{if $i=j$,}\\
0 & \text{otherwise. }
\end{cases}
\end{equation*} 
From Claims 1 and 5, we have
\begin{eqnarray*}
& &(\mbox{the $(i,j)$ entry of }H_1 H_1^T + A'A'^T + B'B'^T) \\
&=&
\begin{cases}
n & \text{if $i-j \in N$,}\\
0 & \text{otherwise,}
\end{cases}
+ 
\begin{cases}
n(n-1) & \text{if $i=j$,}\\
-n & \text{if $i-j \in N$, $i\not=j$,}\\
0 & \text{otherwise,}
\end{cases}
\\
&=& 
\begin{cases}
n^2 & \text{if $i=j$,}\\
0 & \text{otherwise.}
\end{cases}
\end{eqnarray*}
From Claims 1, 4, 5 and 7, we have
\begin{eqnarray*}
& & (\mbox{the $(i,j)$ entry of }H_2 H_2^T + A'A'^T + B'B'^T + 4CC^T+2B'C^T +2C{B'}^T) \\
&=&
\begin{cases}
n & \text{if $i-j \in N$,}\\
0 & \text{otherwise,}
\end{cases}
+
\begin{cases}
n(n-1) & \text{if $i=j$,}\\
-n & \text{if $i-j \in N, i \not= j$,}\\
0 & \text{otherwise,}
\end{cases}\\
& &\hspace{0.4cm}
+ 
\begin{cases}
2n & \text{if $i-j \in N$,}\\
0 & \text{otherwise,}
\end{cases}
+ 
\begin{cases}
-2n &  \text{if $i-j \in N$,}\\
0 & \text{otherwise,}
\end{cases}
\\
&=& 
\begin{cases}
n^2 & \text{if $i=j$,}\\
0 & \text{otherwise.}
\end{cases}
\end{eqnarray*}
We obtain $-J_n H_1^T +  H_1^T B'^T + H_2^T {A'}_1^T = O$, $-H_2J_n + A'H_1 - B'H_2 -2CH_2 = O$
and $H_2 H_1^T + A'B'^T - B'A'^T -2CA'^T = O$ from Claims 8, 9, 10, and $H_1J_n = O$. Finally,  $M$ is symmetric by 
the definitions of $A,B,C$ and  the assumption (i).
Thus we obtain the assertion. 
\qed
%%%%%%%%%
%%%%%%%%%%%%%%%%%%%%%%%%% remark 4.2 %%%%%%%%%%%%%%%%%%%%%%%%%%%%%%
%%%%%%%%%%%%%%%

\vspace{0.2cm}
\begin{remark} 
\begin{itemize}
\item[(1)]
%Goethals and Seidel provided a construction  of symmetric Hadamard matrices  of order $n^2$ under the assumption of the existence of Hadamard matrices of order $n$, see \cite{GS} (or Theorem 5.13 of p.342~\cite{Sze2}). 
%On the other hand, our construction of symmetric Hadamard matrices in Theorem~\ref{Hadamard} assumes the existence of not only Hadamard matrices of order $n$ but also divisible difference families.  
%Hence, our construction will produce symmetric Hadamard matrices nonisomorphic  
%to those obtained by Goethals-Seidel's construction according to the structures of the assumed 
%divisible difference families. 
The known construction theorem by Goethals-Seidel (cf. \cite{GS} or \cite[Theorem~5.13]{Sze2}) provides
the same orders of Hadamard matrices obtained from Theorem~\ref{Hadamard} and requires the existence of an Hadamard matrix only. Though Theorem~\ref{Hadamard} requires further more condition, that is 
the existence of DDFs, we think that Theorem~\ref{Hadamard} may produces Hadamard matrices nonisomorphic to those by Goethals-Seidel's theorem 
according to the assumed DDF. 
\item[(2)] Note that Yamamoto-Yamada's Hadamard difference set
$D=\{a(1+2b)\,|\,a\in \T_n^\ast,\overline{b}\in E\}$ immediately yields 
a symmetric Hadamard matrix of order $2^{2n}$. 
On the other hand, since the divisible difference family obtained from Theorem~\ref{mainthem}
satisfies the conditions of Theorem~\ref{Hadamard}, we have a symmetric Hadamard matrix of order $2^{2n}$. One  can see that these two symmetric Hadamard matrices are equivalent. A similar thing had happened also for the Szekeres difference families: 
if we apply Proposition~\ref{SzeHada} to the Szekeres difference family, we obtain a skew Hadamard matrix of order $q+1$, which is equivalent to that 
obtained from the cyclotomic difference set of index $2$ in $\F_q$ as described in Introduction. 
However, in \cite{Blatt,White}, the authors succeeded to construct difference families satisfying 
the conditions of  Proposition~\ref{SzeHada} different from the Szekeres difference families. Then, they obtained new skew Hadamard matrices by  Proposition~\ref{SzeHada}. Hence, we expect that we have divisible difference families different from those of Theorem~\ref{mainthem} to obtain 
new symmetric Hadamard matrices.  
Thus, we have the following natural problem.  
\begin{problem}
Construct divisible difference families 
satisfying the conditions of Theorem~\ref{Hadamard} different from  divisible 
difference families given in Theorem~\ref{mainthem}. 
\end{problem}
\end{itemize}
\end{remark}
%%%%%%%%%%%%%%%%%%%%%%%%%%%%%%%%%%%%%%%%%%%%%%%%%%%%%%%
\section*{Acknowledgements}
The work of 
K. Momihara was supported by JSPS KAKENHI Grant Number 23840032.
The work of M. Yamada was supported by JSPS KAKENHI Grant Number 90540014.

%%%%%%%%%%%%%%%%%%%%%%%%%%%%%%%%%%%%%%%%%%%%%%%%%%%%%%%


\begin{thebibliography}{25}
%eference
%\bibitem{Beth}
%T. Beth, D. Jungnickel, H. Lenz, {\em Design Theory}, Cambridge University
%Press, (1999).
\bibitem{AHJP} 
K.T. Arasu, W.H. Haemers, D. Jungnickel, A. Pott, Matrix constructions of divisible designs, {\it Lin. Algebra Appl.}, {\bf 153} (1991), 123--133.  
\bibitem{Arasu}
K.T. Arasu, S. Harris, New constructions of group divisible designs, {\it J. Statist. Plann. Inference}, {\bf 52} (1996), 241--253. 
%\bibitem{ADHK}
%K.T. Arasu, C. Ding, T. Helleseth, P.V. Kumar, H.M. Martinsen, 
%Almost difference sets and their sequences with optimal autocorrelation, {\em IEEE Tran. Inform. Theory}, 47, 
%pp. 2934--2943, (2001). 
\bibitem{Blatt} D. Blatt, G. Szekeres, A skew Hadamard matrix of order $52$, {\it Canad. J. Math.}, 
{\bf 21} (1969), 1319--1322.
 
%\bibitem{CD}
%C.J. Colbourn, J.H. Dinitz, {\it Handbook of Combinatorial Designs, 2nd ed.}, 
%Chapman\& Hall/CRC, (2007), G. Ge: {\it Group Divisible Designs,} 255--260.   
\bibitem{Co}
J. Cooper,  A construction for certain classes of supplementary difference sets, {\it J. Austral. Math. Soc.,} {\bf 22} (1976), 182--187.  
%\bibitem{Fer}
%G.A. Fernandez-Alcober, R. Kwashira, L. Martinez, 
%Cyclotomy over products of finite fields and combinatorial applications, {\it Europ. J. Combin.,} 
%{\bf 31} (2010), 1520--1538.
\bibitem{Davis}
J.A. Davis, New semiregular divisible difference sets, {\it Discr. Math.}, {\bf 188} (1998), 99--109. 
\bibitem{GS}
J.-M. Goethals, J.J. Seidel, 
Orthogonal matrices with zero diagonal, {\it Canad. J. Math.}, {\bf 19} (1967),  1001--1010.  
\bibitem{H}
W.H. Haemers, H. Kharaghani, M.A. Meulenberg, Divisible disign graphs, 
{\it J. Combin. Theory, Ser. A}, {\bf 118} (2011), 978--992. 
\bibitem{Hedayat}
A. Hedayat, W.D. Wallis, Hadamard matrices and their applications, 
{\it Ann. Statist.,} {\bf 6} (1978), 1184--1238. 
\bibitem{Horadam}
K.J. Horadam, {\it Hadamard matrices and their applications}, Princeton Univ. Press, (2006). 
\bibitem{Ionin}
Y.J. Ionin, A technique for constructing divisible difference sets, {\it J. Geom.}, {\bf 67} (2000), 164--172.
\bibitem{Kan} L. Kan, {\it Divisible difference sets and difference sets from cyclotomy}, 
Ph.D. Thesis, The Ohio State Univ., (1999). %
\bibitem{J1}
D. Jungnickel, On automorphism groups of divisible designs, {\it Canad. J. Math.,} {\bf 34} (1982), 257--297. 
\bibitem{J2}
D. Jungnickel, On automorphism groups of divisible designs, II: Group invariant generalized conference matrices,  {\it Arch. Math.,} {\bf 54} (1990), 200--208. 

\bibitem{Sarvate}
D.G. Sarvate, J. Seberry, Group divisible designs, GBRSDS and generalized weighing matrices, 
{\it Util. Math.}, {\bf 54} (1998), 157--174.  
\bibitem{Seberry}
J. Seberry, M. Yamada, Hadamard matrices, sequences, block designs, in 
{\em Contemporary design theory: a collection of surveys}, edited by J.H. Dinitz, D.R. Stinson, 
John Wiley \& Sons, Inc, (1992). 
\bibitem{Storer}
T. Storer, {\it Cyclotomy and Difference Sets}, Lectures in Advanced Mathematics, 2, Markham Publishing Company, (1967).
\bibitem{Sze} G. Szekeres, Tournaments and Hadamard matrices, {\it Enseignment Math.}, {\bf 15} (1969), 269--278. 
%
\bibitem{Sze2} W.D. Wallis, A.P. Street, J.S. Wallis, {\it Combinatorics: Room Squares, Sum-Free Sets, 
Hadamard Matrices}, Lecture Notes in Mathematics, 292, Springer-Verlag, (1972).
\bibitem{YY}
K. Yamamoto, M. Yamada,  Hadamard difference sets over an extension of $Z/4Z$, {\it Util. Math.}, {\bf 34} (1988), 169--178. 
\bibitem{WH} 
G. Weng, L. Hu, Some results on skew Hadamard difference sets, {\it Des. Codes and Cryptogr.}, {\bf 50} (2009), 93--105. 
\bibitem{White}
A.L. Whiteman, An infinite family of skew Hadamard matrices, {\it Pacific. J. Math.}, {\bf 38} (1971), 817--822.
\bibitem{XX}
M.Y. Xia, T.B. Xia, Hadamard matrices constructed from supplementary 
difference sets in the class ${\mathcal H}_1$, {\it J. Combin. Des.}, {\bf 2} (1994), 325--339.  
\bibitem{XXS}
M.Y. Xia, T.B. Xia, J. Seberry, A new method for constructing Williamson matrices, {\it Des. Codes, Cryptogr.,} {\bf 35} (2005), 191--209.
\end{thebibliography}
\end{document}